\documentclass[11pt]{article}

\usepackage{float}
\usepackage{fullpage}
\usepackage{amsmath,amssymb,amsthm}
\usepackage{marginnote} %for math mode margin notes.
\usepackage[round]{natbib}
\usepackage{graphicx}
\usepackage{caption}
\usepackage{authblk}

\numberwithin{equation}{section}
\newtheorem{theorem}{Theorem}

\newtheorem{proposition}{Proposition}[section]

\newtheorem{remark}{Remark}[section]

\usepackage{commath}
\usepackage[normalem]{ulem}
\usepackage{mathabx}
\usepackage{mathrsfs}

\usepackage[usenames]{color}
\definecolor{plum}{rgb}{.4,0,.4}
\definecolor{BrickRed}{rgb}{0.6,0,0}
\definecolor{light-gray}{gray}{0.5}
\usepackage[plainpages=false,pdfpagelabels,colorlinks=true,linkcolor=BrickRed,citecolor=plum]{hyperref}

\def\ddefloop#1{\ifx\ddefloop#1\else\ddef{#1}\expandafter\ddefloop\fi}

\def\ddef#1{\expandafter\def\csname c#1\endcsname{\ensuremath{\mathcal{#1}}}}
\ddefloop ABCDEFGHIJKLMNOPQRSTUVWXYZ\ddefloop

\def\ddef#1{\expandafter\def\csname s#1\endcsname{\ensuremath{\mathsf{#1}}}}
\ddefloop ABCDEFGHIJKLMNOPQRSTUVWXYZ\ddefloop

\def\E{\mathop{\mathbf{E}}}

\begin{document}

\title{On the Minimax Optimality of Estimating the Wasserstein Metric}

\author[1]{Tengyuan Liang\thanks{Liang gratefully acknowledges support from the George C. Tiao Fellowship.}}

\affil[1]{University of Chicago, Booth School of Business}

\maketitle

\begin{abstract}
	We study the minimax optimal rate for estimating the Wasserstein-$1$ metric between two unknown probability measures based on $n$ i.i.d. empirical samples from them. We show that estimating the Wasserstein metric itself between probability measures, is not significantly easier than estimating the probability measures under the Wasserstein metric. We prove that the minimax optimal rates for these two problems are multiplicatively equivalent, up to a $\log \log (n)/\log (n)$ factor.
\end{abstract}

\section{Introduction}

In this note we study the minimax optimal rates for estimating the population Wasserstein metric between probability measures based on empirical samples. Let $\mu, \nu$ be two probability measures in $\Omega=[0, 1]^d$, and $W(\mu, \nu)$ denote the Wasserstein-$1$ distance between them. Suppose $X_1, \ldots X_m$ are i.i.d samples from $\mu$, and $Y_1, \ldots, Y_n$ i.i.d from $\nu$. We study:
the minimax optimal rate for estimating $W(\mu, \nu)$ based on $\{X_i\}_{i=1}^m, \{Y_j\}_{j=1}^n$, for some class of probability measures $\cG$ of interest
    \begin{align}
        \inf_{\widetilde T_{m,n}} \sup_{\mu, \nu \in \cG} \E | \widetilde T_{m,n} -  W(\mu, \nu) | \enspace.
    \end{align}
The problem is of importance in both statistics and machine learning, with applications such as nonparametric two sample testing, evaluation of the transportation cost from one set of samples to another, and transfer learning. It turns out that using empirical measures $\widehat \mu_m, \widehat \nu_n$ to estimate is a bad idea. Due to a result by \cite{dudley1969speed}, even for infinitely smooth $\cG = \{{\rm Unif}(\Omega) \}$ and $d\geq 2$,
\begin{align}
     \sup_{\mu, \nu \in \cG} |W(\widehat \mu_m, \widehat\nu_n) -  W(\mu, \nu)| \asymp  n^{-\frac{1}{d}} \enspace.
\end{align}
A natural question arises: can one obtain faster rate, for estimating the Wasserstein metric with other estimators $\widetilde T_{m,n}$ leveraging the regularity of $\cG$ such as smoothness.

A related yet different problem studied in the current literature is estimating a probability measure under the Wasserstein metric based on samples \citep{weed2017sharp, liang2018well, singh2018nonparametric, weed2019estimation}:
    \begin{align}
        \inf_{\widetilde\nu_n} \sup_{\nu \in \cG} \E W(\widetilde\nu_n , \nu) \enspace.
    \end{align}

The two problems are close in nature: ``estimating the metric itself'' is usually an \textbf{easier} problem than ``estimating under the metric.'' In fact, the solution of the latter problem $\widetilde \mu_m, \widetilde \nu_n$ naturally induces a plug-in answer to the first, since
\begin{align*}
    \E |W(\widetilde \mu_m, \widetilde \nu_n) - W(\mu, \nu)| \leq \E W(\widetilde\mu_m , \mu) + \E W(\widetilde\nu_n , \nu) \enspace.
\end{align*}
However, it is unclear whether such a plug-in estimator is optimal. In fact, it is well-known that estimating specific functional of density $F(\nu)$ is usually strictly easier than estimating the density $\nu$ itself. For example, in estimating quadratic functionals of a smooth density vs. estimating under the quadratic functionals, the plug in approach is strictly sub-optimal where the rates can be much improved \citep{bickel1988estimating,donoho1990minimax}.

In this paper, however, we prove that ``estimating the Wasserstein-$1$ metric'', is \textbf{not significantly easier} than ``estimating under the Wasserstein-$1$ metric''. Namely, the plug-in approach is minimax optimal up to a $\log \log (n)/\log (n)$ factor
\begin{align*}
    \frac{\log \log (n \wedge m)}{\log (n \wedge m)} \cdot (n \wedge m)^{-\frac{\beta+1}{2\beta+d}}  \precsim \inf_{\widetilde T_{m,n}} \sup_{\mu, \nu \in \cG_\beta} &\E | \widetilde T_{m,n} -  W(\mu, \nu) |  \\
     \leq  \inf_{\widetilde\mu_m, \widetilde\nu_n} \sup_{\mu, \nu \in \cG_\beta} &\E |W(\widetilde \mu_m, \widetilde \nu_n) - W(\mu, \nu)|  \precsim (n \wedge m)^{-\frac{\beta+1}{2\beta+d}} ,
\end{align*}
where $\cG_\beta$ contains probability measures with densities in H\"{o}lder space with smoothness $\beta \in \mathbb{R}_{\geq 0}$.
The result informs us that seeking other forms of estimators for $W(\mu,\nu)$ would only improve the rates logarithmically. The current result is in contrast with that in a forthcoming companion paper \citep{liang2019adverse}, where we show that ``estimating the adversarial losses'' is \textbf{much easier} than ``estimating under the adversarial losses'', for a collection of integral probability metrics.

Remark that studying the Wasserstein metric and optimal transport for probability measures $\mu, \nu$ with regularity condition has been an important topic in mathematics since Cafferalli's seminal result on regularity theory \citep{caffarelli1991some,caffarelli1992regularity}. By studying the Monge-Amp\'{e}re equation, Cafferalli showed that the Kantorovich potential satisfies specific regularity property, when $\mu, \nu$ are H\"{o}lder smooth.
In this paper, we follow the same H\"{o}lder smooth conditions on $\mu, \nu$, and study the statistical optimal rates for estimating $W(\mu, \nu)$, based on $n$-i.i.d samples.

\subsection{Preliminaries}

Let $\sC^{\beta}(M):=\sC^{\lfloor \beta \rfloor, \beta-\lfloor \beta\rfloor}(M)$ to be H\"{o}lder space with smoothness $\beta \in \mathbb{R}_{\geq 0}$. 
\begin{align}
	\sC^{\beta}(M) :=\left\{ f: \Omega \rightarrow \mathbb{R}: \max_{|\alpha| \leq  \lfloor \beta \rfloor}  \sup_{x \in \Omega} |D^{\alpha} f|  +  \max_{|\alpha| =  \lfloor \beta \rfloor}\sup_{x \neq y \in \Omega}\frac{| D^{\alpha} f(x) - D^{\alpha} f(y) |}{\|x - y \|^{\beta-\lfloor \beta\rfloor}} \leq M \right\}
\end{align}
where $\alpha=[\alpha_1,\ldots,\alpha_d] \in \mathbb{N}^d$ ranges over multi-indices, and $|\alpha| := \sum_{i=1}^d \alpha_i$.
We only consider the bounded case with $\Omega = [0, 1]^d$. 
The class of probability measures of interest is
\begin{align}
	\label{eq:measure-class}
	\cG_\beta:=\left\{ \mu ~:~ \int_\Omega d \mu = 1, \mu \geq 0, \frac{d\mu}{dx} \in \sC^{\beta}(M)   \right\} \enspace.
\end{align}
The Wasserstein-$1$ metric is defined as
\begin{align}
	W_1(\mu, \nu) &:= \inf_{\pi \in \Pi(\mu, \nu)} \int_{\cX \times \cY} \| x - y\| d\pi 
\end{align} 
where $\Pi(\mu, \nu)$ denotes all coupling of probability measures $\mu, \nu$.

\section{Optimal Rates for Estimating Wasserstein Metric}

\begin{theorem}[Minimax Rate]
    Consider $d\geq 2$ and the domain $\Omega = [0, 1]^d$. Given $m$ i.i.d. samples $X_1, \ldots, X_m$ from $\mu$, and $n$ i.i.d. samples from $\nu$, then the minimax optimal rates for estimating $W(\mu,\nu)$ satisfies
    \begin{align}
        \frac{\log \log (n \wedge m)}{\log (n \wedge m)} \cdot (n \wedge m)^{-\frac{\beta+1}{2\beta+d}}  \precsim \inf_{\widetilde T_{m,n}} \sup_{\mu, \nu \in \cG_\beta} &\E | \widetilde T_{m,n} -  W(\mu, \nu) | \precsim (n \wedge m)^{-\frac{\beta+1}{2\beta+d}} \enspace,
    \end{align}
	where the $\mu, \nu$ lies in $\cG_\beta, \beta \geq 0$ as in \eqref{eq:measure-class} whose densities are $\beta$-H\"{o}lder smooth.
\end{theorem}
\begin{remark}
\rm

A few remarks are in order. First, we emphasize that the main technicality is in deriving the lower bound. We construct two composite/fuzzy hypotheses using delicate priors with matching $\log (n \wedge m)$ moments. However, the Wasserstein metric to estimate differs sufficiently under the null vs. under the alternative. Then we calculate the total variation metric directly on the posterior of data defined by the composite hypothesis, using a telescoping technique.

Second, as direct corollary, the following extension hold true. Suppose $\mu \in \cG_{\beta_1}$ and $\nu \in \cG_{\beta_2}$, then define $\beta :=\beta_1 \wedge \beta_2$,
\begin{align}
    \frac{\log \log (n \wedge m)}{\log (n \wedge m)} \cdot (n \wedge m)^{-\frac{\beta+1}{2\beta+d}}  \precsim \inf_{\widetilde T_{m,n}} \sup_{\mu \in \cG_{\beta_1}, \nu \in \cG_{\beta_2}} &\E | \widetilde T_{m,n} -  W(\mu, \nu) | \precsim (n \wedge m)^{-\frac{\beta+1}{2\beta+d}} \enspace.
\end{align}
A further direct implication is: when estimating the cost to transport a known measure $\mu \sim {\rm Unif}([0, 1]^d)$ to an unknown $\nu$ based on $Y_1, \ldots, Y_n$, the result follows from setting $\beta_1 = \infty$ and $m = \infty$.

\end{remark}

\subsection{Proof of the Lower Bound}

Without loss of generality, consider $m\geq n$.
In the lower bound construction, we make use of the multi-resolution analysis. Denote
$\sB^{\beta, p}_{q}$ as the Besov space \citep{tribel1980theory,donoho1996density} with smoothness $\beta \in \mathbb{R}_{\geq 0}$, and $1\leq p, q \leq \infty$,
	\begin{align*}
		\sB^{\beta, p}_{q}(M) := \left\{ f(x) = \sum_{j=0}^{\infty} \sum_{k=0}^{2^{dj}-1}\theta_{jk} h_{jk}(x): \left( \sum_{j=0}^{\infty} \left( (2^{dj})^{s} ( \sum_{k=0}^{2^{dj}-1} |\theta_{jk}|^p )^{1/p} \right)^q  \right)^{1/q} \leq M, ~\text{with}~s = \frac{\beta}{d}+\frac{1}{2}-\frac{1}{p} \right\}
	\end{align*}
	where $h_{jk}(x), x\in [0, 1]^d$ is the wavelet basis.
First, let us review some basic results on function spaces based on
	 \cite{tribel1980theory,donoho1996density}.
	\begin{proposition}
		Under regularity conditions, the following equivalence holds between Besov space and H\"{o}lder space
		\begin{align}
			\sB^{\beta, \infty}_{\infty} = \sC^{\beta}, \text{for}~\beta \notin \mathbb{N}
		\end{align}
		In particular, when $\beta=1$, $\sB^{1,\infty}_{\infty} \supseteq {\rm Lip} \supseteq \sB^{1,\infty}_{1}$.
	\end{proposition}

\paragraph{Step 1: reduction to Besov space norm.}
Write $f_{jk} := \langle f, h_{jk} \rangle$, and $u_{jk} := \langle d \mu/dx, h_{jk}\rangle, v_{jk} := \langle d \nu/dx, v_{jk}\rangle$, we define the following integral probability metric as a surrogate
\begin{align*}
	d_{\sB^{\gamma,p}_q}(\mu, \nu) &:= \sup_{f \in \sB^{\gamma,p}_{q}} |  \int f d\mu  - \int f d\nu |\\
	& = \sup_{f \in \sB^{\gamma,p}_q} | \sum_{j\geq 0} \sum_{k=0}^{2^{dj}-1} f_{jk} (u_{jk} - v_{jk}) | \\
	& = \sup_{f \in \sB^{\gamma,p}_q}  | \sum_{j\geq 0} \| f_{j\cdot} \|_{p} \| u_{j\cdot} - v_{j\cdot} \|_{p_{\star}} | \\
	& = \sup_{f \in \sB^{\gamma,p}_q}  | \sum_{j\geq 0} (2^{dj})^{\frac{\gamma}{d}+\frac{1}{2}-\frac{1}{p}} \| f_{j\cdot} \|_{p} \cdot (2^{-dj})^{\frac{\gamma}{d}+\frac{1}{2}-\frac{1}{p}} \| u_{j\cdot} - v_{j\cdot} \|_{p_{\star}} | \\
	& = \left\{ \sum_{j\geq 0} \left[(2^{dj})^{\frac{\gamma}{d}+\frac{1}{2}-\frac{1}{p}} \| f_{j\cdot} \|_{p} \right]^q \right\}^{1/q} \left\{ \sum_{j\geq 0} \left[(2^{-dj})^{\frac{\gamma}{d}+\frac{1}{2}-\frac{1}{p}} \| u_{j\cdot} - v_{j\cdot} \|_{p_{\star}} \right]^{q_{\star}} \right\}^{1/q_{\star}} \enspace.
\end{align*}
Take $p=q=\infty$ (in this case $p_{\star} = q_{\star} = 1$), we know
\begin{align*}
	d_{\sB^{\gamma,\infty}_\infty}(\mu, \nu) = \sum_{j\geq 0} (2^{-dj})^{\frac{\gamma}{d}+\frac{1}{2}} \sum_{k=0}^{2^j-1} |u_{jk} - v_{jk}|.
\end{align*}
Take $p =\infty$, $q=1$, we know
\begin{align*}
	d_{\sB^{\gamma,\infty}_\infty}(\mu, \nu) = \max_{j \geq 0}(2^{-dj})^{\frac{\gamma}{d}+\frac{1}{2}} \sum_{k=0}^{2^j-1} |u_{jk} - v_{jk}|.
\end{align*}
Now the problem is related to estimation of weighted sum of $\ell_1$ norm of the wavelet coefficients of the densities, in the following multiplicative sense
\begin{align}
  d_{\sB^{\gamma,\infty}_1}(\mu, \nu) \leq 	W(\mu, \nu) \leq d_{\sB^{\gamma,\infty}_\infty}(\mu, \nu) \enspace.
\end{align}
However, multiplicative equivalence is not enough for estimating $W(\mu, \nu)$. In our lower bound construction, we will show that for the hard instances of interest, equality holds.

\paragraph{Step 2: composite hypothesis testing.}
Next we are going to construct two priors on $\nu$ such that 
\begin{align}
	|\E_{\nu \sim \cP_0} W(\mu, \nu) - \E_{\nu \sim \cP_1} W(\mu, \nu)|
\end{align}
are large, while one can not distinguish the following two distributions 
\begin{align}
	p_0(Y_1, \ldots Y_n) =  \E_{\nu \sim \cP_0} \left[ \prod_{i=1}^n \frac{d\nu}{dx}(Y_i) \right],~ p_1(Y_1, \ldots Y_n) =  \E_{\nu \sim \cP_1} \left[ \prod_{i=1}^n \frac{d\nu}{dx}(Y_i) \right]
\end{align}
Here $\cP_0, \cP_1$ are two prior distributions on $\nu$. Consider $\mu$ to be the same distribution under the null $H_0$ and the alternative $H_1$.
Set
\begin{align}
	K \asymp \frac{\log n}{\log \log n}, ~\tau \asymp 1.
\end{align}
The choice will be clear in the later part of the proof. The prior construction is inspired from \cite{lepski1999estimation}, where we borrow the following result.
\begin{proposition}
	\label{prop:hardest-priors}
	For any given positive integer $K$ and $\tau \in \mathbb{R}_{\geq 0}$, there exists two symmetric probability measures $q_0$ and $q_1$ on $[-\tau, \tau]$ such that 
	 \begin{align}
	 	\int_{-\tau}^\tau t^l q_0(dt) = \int_{-\tau}^\tau t^l q_1(dt), ~~ l = 0, 1, \ldots, 2K; \\
		\int_{-\tau}^\tau |t| q_1(dt) - \int_{-\tau}^\tau |t| q_0(dt) = 2\kappa \cdot K^{-1} \tau.
	 \end{align}
	 where $\kappa$ is some constant depending on $K$ only.
\end{proposition}

Now let's construct $\cP_0$ and $\cP_1$ as follows. Take $\mu \sim {\rm Unif}([0,1]^d)$. Choose $J \in \mathbb{N}_{\geq 0}$ such that $2^{dJ} \asymp n^{\frac{1}{1+2\beta/d}}$, first we are going to embed a parametrized class of densities into $\sC^{\beta}$
\begin{align}
	\frac{d\nu_{\theta}}{dx} := \mu(x) + \frac{1}{\sqrt{n}} \sum_{k=0}^{2^{dJ}-1} \theta_k h_{Jk}(x)
\end{align}
with $\theta_{k} \in [-\tau, \tau]$ for all $k$.

We will now show that the construction lies inside the measure class $\nu_\theta \in \cG_{\beta}$.
First observe that for wavelet basis that satisfy the regularity condition $\int h_{jk} d\mu = 0$, we have $\int_{\Omega} \nu_\theta dx = 1$ and $d \nu_\theta/dx \geq 1 - \sqrt{2^{dJ}/n} >0$. Hence it is a valid probability measure.
Let's then verify $d \nu_\theta/ dx \in \sB^{\beta, \infty}_{1} \subseteq \sC^{\beta}$ lies in the H\"{o}lder space. This follows since
\begin{align}
	\frac{1}{\sqrt{n}} |\theta_k| \leq (2^{dJ})^{-(\frac{\beta}{d}+\frac{1}{2})}, ~~\forall k.
\end{align} 

For any $\gamma \geq 0$
\begin{align*}
	d_{\sB^{\gamma,\infty}_\infty}(\mu, \nu_{\theta})&:=  (2^{-dJ})^{\frac{\gamma}{d}+\frac{1}{2}} \frac{1}{\sqrt{n}}\sum_{k=0}^{2^{dJ}-1} |\theta_k| \\
	&=(2^{-dJ})^{\frac{\gamma}{d}+\frac{1}{2}} (2^{dJ})^{-(\frac{\beta}{d}+\frac{1}{2})} \sum_{k=0}^{2^{dJ}-1} |\theta_k| \\
	&=(2^{-dJ})^{\frac{\beta+\gamma}{d}} \frac{1}{2^{dJ}}\sum_{k=0}^{2^{dJ}-1} |\theta_k|.
\end{align*}
It is easy to verify that
\begin{align*}
	 d_{\sB^{\gamma,\infty}_1}(\mu, \nu_{\theta}) = (2^{-dJ})^{-\frac{\beta+\gamma}{d}} \cdot \frac{1}{2^{dJ}} \sum_{k\in [2^{dJ}]}|\theta_k| = d_{\sB^{\gamma,\infty}_\infty}(\mu, \nu_{\theta}) 
\end{align*}
Therefore we must have for any $q \geq 1$, take $\gamma = 1$
\begin{align*}
	W(\mu, \nu_\theta) =  d_{\sB^{1,\infty}_q}(\mu, \nu_{\theta}) = (2^{-dJ})^{-\frac{\beta+1}{d}} \cdot \frac{1}{2^{dJ}} \sum_{k\in [2^{dJ}]}|\theta_k| .
\end{align*}

\paragraph{Step 3: polynomials and matching moments.}
Recall the collection of measures $\cS_{0} := \{\nu_\theta: \theta_k \sim q_0~ i.i.d.~ \text{for}~ k\in [2^{dJ}] \}$, and $\cP_0$ can be viewed as an uniform prior over this set $\cS_{0}$. Similar construction for $\cP_1$ via $q_1$.
Remark that due to the separation of support for wavelets, we have
\begin{align}
	\frac{d\nu_{\theta}}{dx} = \prod_{k=1}^{2^{dJ}} (1+ \theta_k n^{-1/2} h_{Jk}(x)) \enspace.
\end{align}
Therefore we know
\begin{align}
	\label{eq:sum-prod}
	p_0(Y_1,\ldots,Y_n) &= \E_{\theta \sim q_0^{\otimes 2^{dJ}}}  \prod_{i=1}^n \frac{d\nu_{\theta}}{dx}(Y_i) = \E_{\theta \sim q_0^{\otimes 2^{dJ}}}  \prod_{i=1}^n \prod_{k=1}^{2^{dJ}} (1+\theta_k n^{-1/2} h_{Jk}(Y_i)) \\
	&= \E_{\theta \sim q_0^{\otimes 2^{dJ}}}  \prod_{k=1}^{2^{dJ}} \prod_{i=1}^n (1+\theta_k n^{-1/2} h_{Jk}(Y_i)) \\
	&= \prod_{k=1}^{2^{dJ}} \E_{\theta_k \sim q_0} \prod_{i=1}^n (1+\theta_k n^{-1/2} h_{Jk}(Y_i)) \enspace.
\end{align}

Let's analyze the polynomial in $\theta_k$ (and $h_{Jk}(Y_i)$) with degree at most $n$
\begin{align}
	f(\theta_k;h_{jk}(Y_1),\ldots,h_{jk}(Y_n))&:=  \prod_{i=1}^n (1+\theta_k \frac{h_{Jk}(Y_i)}{\sqrt{n}}) \\
	&=\sum_{l=0}^{n} \theta_k^l \frac{\sum_{i_1<\ldots<i_l} h_{Jk}(Y_{i_1})\ldots h_{Jk}(Y_{i_l})}{n^{l/2}}\\
	&=: \sum_{l=0}^{n} \theta_k^l \frac{H^{(l)}_{Jk}(Y_1, \ldots, Y_n)}{n^{l/2}} \label{eq:polyn}
\end{align}
where $H^{(l)}_{JK}(Y_1, \ldots, Y_n)$ a sum of monomial of order $l$, i.e., $\binom{n}{l}$ terms with each of the form $h_{Jk}(Y_{i_1})\ldots h_{Jk}(Y_{i_l})$.
Denote $f^{[\leq K]}, f^{[> K]}$ to denote the corresponding truncated polynomial according to degree.

In this convenient notation, we know
\begin{align}
	p_0(Y_1,\ldots,Y_n) = \prod_{k \in [2^{dJ}]} \E_{\theta_k \sim q_0}f(\theta_k;h_{Jk}(Y_1),\ldots,h_{Jk}(Y_n))
\end{align}
Later, we shall use the following properties of the polynomial $f$ of degree at most $n$.
\begin{align}
	\forall \theta_k, ~~\int_{\cY^{\otimes n}} f(\theta_k;h_{Jk}(y_1),\ldots,h_{Jk}(y_n)) dy_1\ldots dy_n= 1 % \\
% 	\forall y_1,\ldots,y_n, ~~\int_{[-\tau,\tau]} f(\theta;h_{jk}(Y_1),\ldots,h_{jk}(Y_n)) q_0(d\theta)
\end{align}
And the following property according to $q_0$ and $q_1$ constructed in Proposition~\ref{prop:hardest-priors}: $\forall y_1,\ldots,y_n$
\begin{align*}
	&\E_{\theta_k \sim q_1} f(\theta_k;h_{Jk}(y_1),\ldots,h_{Jk}(y_n)) - \E_{\theta_k \sim q_0} f(\theta_k;h_{Jk}(y_1),\ldots,h_{Jk}(y_n)) \\
	&= \int_{[-\tau,\tau]} f^{[>2K]}(\theta_k;h_{Jk}(y_1),\ldots,h_{Jk}(y_n)) (q_1-q_0)(d\theta_k) \enspace.
\end{align*}

\paragraph{Step 4: total variation and telescoping.}

\begin{align*}
	{\rm TV}(p_1, p_0)&:= \frac{1}{2} \int_{\cY^{\otimes n}} \left|p_1(y_1, \ldots, y_n) - p_0(y_1, \ldots, y_n) \right| dy_1\ldots dy_n \\
	&= \frac{1}{2} \int_{\cY^{\otimes n}}  \left|   \prod_{k\in [2^{dJ}]} E_{\theta_k \sim q_1} f(\theta_k;h_{Jk}(y^{\otimes n})) - \prod_{k \in [2^{dJ}]} \E_{\theta_k \sim q_0}f(\theta_k;h_{Jk}(y^{\otimes n})) \right|  dy_1\ldots dy_n \\
\end{align*}

Claim the following telescoping lemma holds. The proof can be done through induction.
\begin{proposition}
	For all $a_i, b_i \geq 0$,
	\begin{align}
		|\prod_{k\in[1, N]} a_k - \prod_{k \in [1, N]} b_k| \leq \sum_{i\in [1, N]} |a_i - b_i| \cdot \prod_{k\in [1, i)} b_k  \cdot \prod_{k\in (i, N]} a_k \enspace.
	\end{align}
\end{proposition}

Define
\begin{align}
	a_k(h_{Jk}(y_1),\ldots, h_{Jk}(y_n)) := E_{\theta_k \sim q_1} f(\theta_k;h_{Jk}(y^{\otimes n})) \\
	b_k(h_{Jk}(y_1),\ldots, h_{Jk}(y_n)) := E_{\theta_k \sim q_0} f(\theta_k;h_{Jk}(y^{\otimes n}))
\end{align}
Using the the above telescoping proposition, we have
\begin{align}
	{\rm TV}(p_1, p_0)&\leq \sum_{k \in [2^{dJ}]} \int |a_k - b_k| \cdot \prod_{k'\in [1, k)} b_{k'}  \prod_{k''\in (k, N]} a_{k''} dy^{\otimes n} \\
	&= \sum_{k \in [2^{dJ}]} \E_{\substack{\theta_{k'} \sim q_0, k' \in [1, k)\\
	\theta_{k''} \sim q_1, k'' \in (k, 2^{dJ}]} } \E_{Y_1,\ldots,Y_n \sim \nu_{\theta_{-k}}} |a_k(h_{Jk}(Y_1),\ldots, h_{Jk}(Y_n)) - b_k(h_{Jk}(Y_1),\ldots, h_{Jk}(Y_n))|
\end{align}
Let's analyze the term
\begin{align*}
	\E_{Y_1,\ldots,Y_n \sim \nu_{\theta_{-k}}} |a_k(h_{Jk}(Y_1),\ldots, h_{Jk}(Y_n)) - b_k(h_{Jk}(Y_1),\ldots, h_{Jk}(Y_n))|
\end{align*}
where $Y_1, \ldots Y_n$ i.i.d. sampled from a measure
\begin{align}
	d \nu_{\theta_{-k}}/dx := 1+ \frac{1}{\sqrt{n}} \sum_{k'\neq k} \theta_{k'} h_{J k'}(x) \enspace.
\end{align}
Note that $\nu_{\theta_{-k}}$ agrees with the uniform measure $\mu$ on the domain associated with $h_{Jk}(x)$.
Due to the separation of support for wavelet basis, we know the random variables
\begin{align}
	h_{J k}(Y_i)
\end{align}
are only determined by $\nu_{\theta_{-k}}$ restricted to the domain of $h_{Jk}$. Hence for $Y_1, \ldots, Y_n \sim \nu_{\theta_{-k}}$, 
\begin{align*}
	&\E_{Y_1,\ldots,Y_n \sim \nu_{\theta_{-k}}} |a_k(h_{Jk}(Y_1),\ldots, h_{Jk}(Y_n)) - b_k(h_{Jk}(Y_1),\ldots, h_{Jk}(Y_n))|\\
	& = \E_{Y_1,\ldots,Y_n \sim \mu} |a_k(h_{Jk}(Y_1),\ldots, h_{Jk}(Y_n)) - b_k(h_{Jk}(Y_1),\ldots, h_{Jk}(Y_n))| \enspace.
\end{align*}
Now one can directly bound the TV metric between the complex sum-product distribution $p_0$ and $p_1$ defined in \eqref{eq:sum-prod},
\begin{align}
	\label{eq:teles-later}
	2{\rm TV}(p_1, p_0)& \leq \sum_{k=1}^{2^{dJ}} \E_{Y_1,\ldots,Y_n \sim \mu} |a_k(h_{Jk}(Y_1),\ldots, h_{Jk}(Y_n)) - b_k(h_{Jk}(Y_1),\ldots, h_{Jk}(Y_n))| \\
	& = \sum_{k=1}^{2^{dJ}} \int \left| \E_{\theta_k \sim q_1} f(\theta_k;h_{Jk}(y^{\otimes n})) - \E_{\theta_k \sim q_0}f(\theta_k;h_{Jk}(y^{\otimes n}))\right| dy_1 \ldots dy_n .
\end{align}

\paragraph{Step 5: $\ell_2$ bound.}
In this section, we are going to bound, for a fixed $k$, the following expression using the properties of the $q_1$ and $q_0$ constructed with matching moments up to $2K$,
\begin{align*}
	\int \left| \E_{\theta_k \sim q_1} f(\theta_k;h_{Jk}(y^{\otimes n})) - \E_{\theta_k \sim q_0}f(\theta_k;h_{Jk}(y^{\otimes n}))\right| dy_1 \ldots dy_n \enspace.
\end{align*}
First, observe the $\ell_2$ bound
\begin{align}
	\int |g_1 -g_2| d\mu \leq \left( \int (g_1 - g_2)^2 d\mu \right)^{1/2} 
\end{align}
Let's bound the $\ell_2$ form
\begin{align}
	&\int \left( \E_{\theta_k \sim q_1} f(\theta_k;h_{Jk}(y^{\otimes n})) - \E_{\theta_k \sim q_0}f(\theta_k;h_{Jk}(y^{\otimes n}))\right)^2 dy_1 \ldots dy_n \label{eq:l2}\\
	&= \E_{\theta, \theta'\sim q_1} \int f(\theta;h_{Jk}(y^{\otimes n})) f(\theta';h_{Jk}(y^{\otimes n})) d y^{\otimes n} + \E_{\omega, \omega'\sim q_0} \int f(\omega;h_{Jk}(y^{\otimes n})) f(\omega';h_{Jk}(y^{\otimes n})) d y^{\otimes n} \nonumber\\
	& \quad - 2 \E_{\theta \sim q_1, \omega\sim q_0}\int f(\theta;h_{Jk}(y^{\otimes n})) f(\omega;h_{Jk}(y^{\otimes n})) d y^{\otimes n} \nonumber
\end{align}
Note now each $f(\theta_k;h_{Jk}(y^{\otimes n}))f(\theta';h_{Jk}(y^{\otimes n}))$ for fixed $\theta, \theta'$ takes the following product form
\begin{align*}
	f(\theta_k;h_{Jk}(y^{\otimes n}))f(\theta';h_{Jk}(y^{\otimes n})) = \prod_{i=1}^n \left(1+ (\theta + \theta') \frac{h_{Jk}(Y_i)}{\sqrt{n}} + \theta \theta' \frac{h^2_{Jk}(Y_i)}{n}  \right)
\end{align*}
and
\begin{align*}
	\int f(\theta;h_{Jk}(y^{\otimes n})) f(\theta';h_{Jk}(y^{\otimes n})) d y^{\otimes n} &= \left( 1 + \theta\theta' \frac{\int h^2_{Jk}(y) dy}{n} \right)^n \\
	&= \left( 1 + \theta\theta' \frac{1}{n} \right)^n \enspace.
\end{align*}

Therefore we have for \eqref{eq:l2} 
\begin{align*}
	\eqref{eq:l2}&= \E_{\theta, \theta'\sim q_1} \left[\left( 1 + \theta\theta' \frac{1}{n} \right)^n \right] + \E_{\omega, \omega'\sim q_0} \left[\left( 1 + \omega\omega' \frac{1}{n} \right)^n \right] -2 \E_{\theta \sim q_1, \omega\sim q_0} \left[\left( 1 + \theta\omega \frac{1}{n} \right)^n \right] \\
	&= \sum_{l=1}^{\lfloor n/2\rfloor}  \left(\E_{\theta, \theta'\sim q_1} [(\theta \theta')^{2l}] + \E_{\omega, \omega'\sim q_0} [(\omega \omega')^{2l}] - 2 \E_{\theta \sim q_1, \omega\sim q_0} [(\theta \omega)^{2l}] \right) \frac{\binom{n}{2l}}{n^{2l}} \\
	&= \sum_{l=1}^{\lfloor n/2\rfloor}  \left( \left(\E_{q_1} [\theta^{2l}] \right)^2 + \left(\E_{q_0} [\theta^{2l}] \right)^2 - 2 \E_{q_1} [\theta^{2l}]  \E_{q_0} [\theta^{2l}] \right) \frac{\binom{n}{2l}}{n^{2l}} 
\end{align*}
Recall the crucial property that for all $l\leq K$, we know
\begin{align}
	\E_{\theta \sim q_1} [\theta^{2l}] = \E_{\theta \sim q_0} [\theta^{2l}] ~~ \Rightarrow~~
	\left(\E_{q_1} [\theta^{2l}] \right)^2 + \left(\E_{q_0} [\theta^{2l}] \right)^2 - 2 \E_{q_1} [\theta^{2l}]  \E_{q_0} [\theta^{2l}] = 0
\end{align}
therefore the above summation equals
\begin{align*}
	\eqref{eq:l2}&=\sum_{l=K+1}^{\lfloor n/2\rfloor} \left( \left(\E_{q_1} [\theta^{2l}] \right)^2 + \left(\E_{q_0} [\theta^{2l}] \right)^2 - 2 \E_{q_1} [\theta^{2l}]  \E_{q_0} [\theta^{2l}] \right) \frac{\binom{n}{2l}}{n^{2l}} \\
	&\leq \sum_{l=K+1}^{\lfloor n/2\rfloor}  4\tau^{4l} \frac{1}{(2l)!} \\
	& \precsim 4\frac{\tau^{4K}}{(2K)!} \exp(\tau^4) \enspace.
\end{align*}

Assemble the two bounds, we have
\begin{align}
	\label{eq:crucial-ineq}
	&\int \left| \E_{\theta_k \sim q_1} f(\theta_k;h_{Jk}(y^{\otimes n})) - \E_{\theta_k \sim q_0}f(\theta_k;h_{Jk}(y^{\otimes n}))\right| dy_1 \ldots dy_n \\
	&\leq 2 \frac{\tau^{2K}}{\sqrt{(2K)!}} \exp(\tau^4/2)
\end{align}

\paragraph{Step 6: combine all pieces.}

Now continuing \eqref{eq:teles-later}, we have
\begin{align}
	2{\rm TV}(p_1, p_0)& \leq \sum_{k=1}^{2^{dJ}} \E_{Y_1,\ldots,Y_n \sim \mu} |a_k(h_{Jk}(Y_1),\ldots, h_{Jk}(Y_n)) - b_k(h_{Jk}(Y_1),\ldots, h_{Jk}(Y_n))| \\
	& = \sum_{k=1}^{2^{dJ}} \int \left| \E_{\theta_k \sim q_1} f(\theta_k;h_{Jk}(y^{\otimes n})) - \E_{\theta_k \sim q_0}f(\theta_k;h_{Jk}(y^{\otimes n}))\right| dy_1 \ldots dy_n \\
	&\leq 2^{dJ}  \cdot 2 \frac{\tau^{2K}}{\sqrt{2K}!} \exp(\tau^4/2) \precsim \exp (c \log n - K \log K)
\end{align}
Therefore by taking $K = \frac{c}{2} \frac{\log n}{\log \log n}$, we know
\begin{align}
	2{\rm TV}(p_1, p_0) \leq n^{-\frac{c}{2}\log n} \leq n^{-c/2}.
\end{align}

We know by construction of the composite hypothesis 
\begin{align*}
	& \quad |\E_{\nu_{\theta} \sim \cP_0} d_{\sB^{\gamma,\infty}_q}(\mu, \nu_{\theta}) - \E_{\nu_{\theta} \sim \cP_1} d_{\sB^{\gamma,\infty}_q}(\mu, \nu_{\theta})| \\
	&=  (2^{-dJ})^{-\frac{\beta+\gamma}{d}} \cdot \left| \E_{\nu_{\theta} \sim \cP_0} \left[ \frac{1}{2^{dJ}} \sum_{k\in [2^{dJ}]}|\theta_k| \right] - \E_{\nu_{\theta} \sim \cP_1} \left[ \frac{1}{2^{dJ}} \sum_{k\in [2^{dJ}]}|\theta_k| \right]   \right|  \\
	& = n^{-\frac{\beta+\gamma}{2\beta+d}} \cdot \left| \E_{\theta \sim q_0} [|\theta|] -\E_{\theta \sim q_1}[ |\theta|] \right| \\
	& \succsim n^{-\frac{\beta+\gamma}{2\beta+d}} \cdot 2\kappa K^{-1} \tau  = n^{-\frac{\beta+\gamma}{2\beta+d}} \cdot \frac{\log \log (n)}{\log(n)} \enspace.
\end{align*}

Therefore we have for any functional of $\theta$, for any estimator based on $n$-i.i.d. samples
\begin{align*}
	\sup_{\nu_\theta}\E_{\cD_n \sim \theta} |\hat{T}_n - F(\theta)| &\geq \E_{\theta \sim Q_0} \E |\hat{T}_n - F(\theta)| \\
	&\geq \E_{\theta \sim Q_0} \E_{\cD_n \sim \theta} |\hat{T}_n - \E_{\theta\sim Q_0}F(\theta)| - \delta_{Q_0}
\end{align*}
where $\delta_{Q_0}:=\E_{\theta \sim Q_0} |\E_{\theta\sim Q_0}F(\theta) - F_{\theta}| $. Here $Q_0$ is some prior distribution on $\theta$.
Repeat the same argument for $Q_1$,
and by Le Cam's argument on two composite hypothesis
\begin{align*}
	\sup_{\nu_\theta} \E |\hat{T}_n - F(\theta)|& \geq \frac{1}{2}\left( \E_{\theta \sim Q_0} \E_{\cD_n \sim \theta} |\hat{T}_n - \E_{\theta\sim Q_0}F(\theta)| + \E_{\theta \sim Q_1} \E_{\cD_n \sim \theta} |\hat{T}_n - \E_{\theta\sim Q_1}F(\theta)| \right) - \frac{\delta_{Q_0} + \delta_{Q_1}}{2} \\
	& = \frac{1}{2}\left( \E_{\cD_n \sim p_0}|\hat{T}_n - \E_{\theta\sim Q_0}F(\theta)| + \E_{\cD_n \sim p_1} |\hat{T}_n - \E_{\theta\sim Q_1}F(\theta)| \right) - \frac{\delta_{Q_0} + \delta_{Q_1}}{2} \\
	& \geq \frac{|\E_{\theta\sim Q_0}F(\theta) - \E_{\theta\sim Q_1}F(\theta)|}{4} \left( P_{0}(T = 1) + P_{1}(T = 0)\right) - \frac{\delta_{Q_0} + \delta_{Q_1}}{2} \\
	& \geq \frac{|\E_{\theta\sim Q_0}F(\theta) - \E_{\theta\sim Q_1}F(\theta)|}{4} \int p_0(y^{\otimes n}) \wedge p_1(y^{\otimes n}) dy^{\otimes n} -  \frac{\delta_{Q_0} + \delta_{Q_1}}{2} \\
	& = \frac{|\E_{\theta\sim Q_0}F(\theta) - \E_{\theta\sim Q_1}F(\theta)|}{4} (1 - d_{TV}(p_0, p_1)) -    \frac{\delta_{Q_0} + \delta_{Q_1}}{2} 
\end{align*}
where $p_i(y^{\otimes n}) = \int Pr(y^{\otimes n} |\theta) Q_i(d\theta)$, for $i=0,1$. Here the test $T=1$ if and only if $\hat{T}_n$ is closer to $\E_{\theta\sim Q_1}F(\theta)$.
In our case, for any $q\geq 1$
$$
F(\theta) :=W(\mu, \nu) = d_{\sB^{1,\infty}_q}(\mu, \nu_{\theta}) = (2^{-dJ})^{-\frac{\beta+1}{d}} \left[ \frac{1}{2^{dJ}} \sum_{k\in [2^{dJ}]}|\theta_k| \right]
$$
then
\begin{align*}
	|\E_{\theta\sim Q_0}F(\theta) - \E_{\theta\sim Q_1}F(\theta)| &= |\E_{\nu_{\theta} \sim \cP_0} d_{\sB^{1,\infty}_q}(\mu, \nu_{\theta}) - \E_{\nu_{\theta} \sim \cP_1} d_{\sB^{1,\infty}_q}(\mu, \nu_{\theta})| \\
	&\succsim n^{-\frac{\beta+1}{2\beta+d}} \cdot \frac{\log \log (n)}{\log(n)} \\
	 1 - d_{TV}(p_0, p_1) &\geq 1 -  n^{-c/2}\\
	 \frac{\delta_{Q_0} + \delta_{Q_1}}{2} &\precsim n^{-\frac{\beta+1}{2\beta+d}} \frac{1}{\sqrt{2^{dJ}}} \ll n^{-\frac{\beta+1}{2\beta+d}} \cdot \frac{\log \log (n)}{\log(n)} \enspace.
\end{align*}
Therefore we have
\begin{align}
		\inf_{\widehat{T}_n} \sup_{\nu \in \sC^{\beta}} \E|\widehat{T}_n - W(\mu, \nu)| \succsim  n^{-\frac{\beta+1}{2\beta+d}} \cdot \frac{\log \log (n)}{\log(n)} \enspace.
\end{align}

% \end{proof}

\subsection{Proof of the Upper Bound}

The upper bound can be obtained through similar derivations as in \cite{liang2018well, singh2018nonparametric, weed2019estimation}. We include here for completeness. 

The estimator is of the plug-in form, with 
\begin{align}
	W(\widetilde \mu_m, \widetilde \nu_n) := \sup_{f \in {\rm Lip(1)}} |\int f d \widetilde{\mu}_m - \int f d\widetilde{\nu}_n |
\end{align}
where $\widetilde \mu_m$, and $\widetilde \nu_n$ are smoothed empirical measures based on truncation on Wavelets.
It is clear that
\begin{align}
	|W(\widetilde \mu_m, \widetilde \nu_n) - W(\mu, \nu)| \leq \sup_{f \in {\rm Lip(1)}} |\int f d \widetilde{\mu}_m - \int f d\mu | + \sup_{f \in {\rm Lip(1)}} |\int f d \widetilde{\nu}_n - \int f d\nu |.
\end{align}

Now let's bound $\sup_{f \in {\rm Lip(1)}} |\int f d \widetilde{\nu}_n - \int f d\nu |$ via expanding under the Wavelet basis. Denote $\widehat{\E}[h_{jk}]:=1/n\sum_{i=1}^n h_{jk}(Y_i)$,  the smoothed empirical estimate $\widetilde\nu_n$ is defined
\begin{align}
	\frac{d \widetilde\nu_n}{dx} := \sum_{j=0}^{J} \sum_{k=0}^{2^{dj}-1} \widehat{\E}[h_{jk}] h_{jk}(x) \enspace.
\end{align}
Expand $f(x) = \sum_{j\geq 0} \sum_{k=0}^{2^{dj}-1} f_{jk} h_{jk}(x)$, we have
\begin{align*}
	 & \sup_{f \in {\rm Lip(1)}} |\int f d \widetilde{\nu}_n - \int f d\nu | \leq \sup_{f \in {\rm \sB^{1,\infty}_{\infty}}} |\int f d \widetilde{\nu}_n - \int f d\nu | \\
	 & = \sup_{f \in {\rm \sB^{1,\infty}_{\infty}}} | \sum_{j \geq 0}^{J}  \sum_{k=0}^{2^{dj}-1} f_{jk} (\widehat{\E}[h_{jk}] - \E[h_{jk}])| +  \sup_{f \in {\rm \sB^{1,\infty}_{\infty}}} | \sum_{j >J}  \sum_{k=0}^{2^{dj}-1} f_{jk}  \E[h_{jk}] |
\end{align*}
For the first term, since $f \in \sB^{1, \infty}_{\infty} \Rightarrow \forall j,k, ~|f_{jk}| \leq (2^{-dj})^{\frac{1}{d}+\frac{1}{2}}$
\begin{align*}
	&\E \sup_{f \in {\rm \sB^{1,\infty}_{\infty}}} | \sum_{j \geq 0}^{J}  \sum_{k=0}^{2^{dj}-1} f_{jk} (\widehat{\E}[h_{jk}] - \E[h_{jk}])| \leq \sum_{j\geq 0}^J (2^{-dj})^{\frac{1}{d}+\frac{1}{2}} \sum_{k=0}^{2^{dj}-1} \E|\widehat{\E}[h_{jk}] - \E[h_{jk}]|\\
	&\leq \sum_{j\geq 0}^J (2^{-dj})^{\frac{1}{d}+\frac{1}{2}} \sum_{k=0}^{2^{dj}-1} (\E|\widehat{\E}[h_{jk}] - \E[h_{jk}]|^2)^{1/2}\\
	&\precsim \sum_{j\geq 0}^J (2^{-dj})^{\frac{1}{d}+\frac{1}{2}} 2^{dj} \frac{1}{\sqrt{n}} \asymp \frac{1}{\sqrt{n}}(2^{dJ})^{\frac{1}{2} - \frac{1}{d}} 
\end{align*}
for $d\geq 2$.

For the second term, recall $\E_{Y\sim \nu} [h_{jk}(Y)] = \langle d\nu/dx, h_jk \rangle =: v_{jk}$. Due to the fact that
\begin{align}
	d \nu/dx \in \sC^{\beta} \in \sB^{\beta, \infty}_{\infty}  \Rightarrow \forall j,k, ~|v_{jk}| \leq (2^{-dj})^{\frac{\beta}{d}+\frac{1}{2}}\\
	f \in \sB^{1, \infty}_{\infty} \Rightarrow \forall j,k, ~|f_{jk}| \leq (2^{-dj})^{\frac{1}{d}+\frac{1}{2}}
\end{align}
\begin{align*}
	& \E \sup_{f \in {\rm \sB^{1,\infty}_{\infty}}} | \sum_{j >J}  \sum_{k=0}^{2^{dj}-1} f_{jk} \E[h_{jk}] | = \E \sup_{f \in {\rm \sB^{1,\infty}_{\infty}}} | \sum_{j >J}  \sum_{k=0}^{2^{dj}-1} f_{jk} v_{jk} | \\
	& \leq  \sum_{j >J}  \sum_{k=0}^{2^{dj}-1} (2^{-dj})^{\frac{1}{d}+\frac{1}{2}} (2^{-dj})^{\frac{\beta}{d}+\frac{1}{2}} \\
	& \leq (2^{dJ})^{-\frac{\beta+1}{d}}
\end{align*}

Balancing the two terms, we have
\begin{align}
	& \sup_{\nu \in \cG_\beta} \sup_{f \in {\rm Lip(1)}} |\int f d \widetilde{\nu}_n - \int f d\nu | \precsim \frac{1}{\sqrt{n}}(2^{dJ})^{\frac{1}{2} - \frac{1}{d}}  + (2^{dJ})^{-\frac{\beta+1}{d}}\\
	& \asymp n^{-\frac{\beta+1}{2\beta+d}}, ~~\text{with $2^{dJ} \asymp n^{\frac{1}{2\beta/d+1}}$}\enspace.
\end{align}
Put everything together, we know
\begin{align}
	\E|W(\widetilde \mu_m, \widetilde \nu_n) - W(\mu, \nu)| \leq (n \wedge m)^{-\frac{\beta+1}{2\beta+d}}.
\end{align}

\bibliographystyle{Chicago}
\bibliography{ref}
\nocite*{}

\end{document}